\documentclass[12pt]{amsart}

\usepackage{amssymb}

\newcommand{\p}{{\mathfrak p}}
\newcommand{\oo}{{\mathcal O}}
\newcommand{\rarr}{{\rightarrow}}

\newtheorem*{lem}{Lemma}

\advance \topmargin by -\headheight
\advance \topmargin by -\headsep
     
\evensidemargin \oddsidemargin
\begin{document}

\title{Fiber products of hyperelliptic curves}

\author{Siman Wong}

\address{Department of Mathematics \& Statistics, University of Massachusetts.
	Amherst, MA 01003-4515 USA}

\email{siman@math.umass.edu}

\maketitle

\thispagestyle{empty}

Let $k$ be a number field, and let $S$ be a finite set of maximal ideals of
the ring of integers $\oo_k$ of $k$.  In his 1962 ICM address, Shafarevich
\cite{shaf} asked if there are only finitely many $k$-isomorphism classes of
algebraic curves of a fixed genus $g\ge 1$ with good reduction outside $S$.
He verified this for $g=1$ by reducing the problem to
Siegel's theorem.  Parshin  \cite{parshin} extended this argument to all
hyperelliptic curves (cf.~also \cite{oort}).  The general case was settled by
Faltings' celebrated work \cite{faltings}.  In this note we give a short proof
of Shafarevich's
conjecture for hyperelliptic curves, by reducing the problem to the case
$g=1$ using the Theorem of de Franchis plus standard facts about
discriminants of hyperelliptic equations.

Let $k, \oo_k$ and $S$ be as above. Enlarge $S$ by a finite number of maximal
ideals if necessary, so that
\newcounter{property}
\begin{list}
	{(\roman{property})}{\usecounter{property}
				\setlength{\labelwidth}{15pt}}
\item
$S$ contains every maximal ideal of $\oo_k$ lying above $2$; and
\item
the ring $\oo_S$ of $S$-integers in $k$ is a principal ring.
\end{list}
Denote by $k_\p$ the residue field of a maximal ideal $\p$ of $\oo_S$. Denote
by $\overline{F}$ a separable closure of a field $F$.

Fix an integer $g\ge 2$.
Denote by $k(g, S)$ the compositum of all extensions of $k$ of degree
dividing $(2g+2)!$ which are unramified outside $S$.  By Hermite's theorem,
$
k(g, S)/k
$
is a finite extension.  To prove Shafarevich's Conjecture for hyperelliptic
curves we proceed by induction on $g$,
assuming as known the case $g=1$.  Following Lockhart \cite{lock}, we say that
a hyperelliptic curve over
$k$ is \textit{pointed} if it has a $k$-rational Weierstrass point.

Denote by $\Sigma(k, g, S)$ the set of hyperelliptic curves over $k$ of genus
$g$ with good reduction outside $S$.
Define $\Sigma_0(k, g, S) = \{ X\in\Sigma(k, g, S): \text{$X$ is pointed over
$k$} \}$.

\begin{lem}
Let $g\ge 2$. To prove that $\Sigma(k, g, S)$ is finite for every $k$
and $S$,
it suffices to prove that $\Sigma_0(k, g, S)$ is finite for every $k$ and
$S$.

\end{lem}

\begin{proof}
If $X/k$ is a curve of genus $\ge 2$, then $\text{Aut}_{k(g, S)}(X)$ is finite,
and hence so does
$
H^1(\text{Gal}(k(g, S)/k), \text{Aut}_{k(g, S)}(X) )
$.
Thus there are at most finitely many $k$-isomorphism classes of
curves over $k$ which become isomorphic to $X$ over $k(g, S)$.  To prove the
lemma it then suffices to show that the Weierstrass points of every curve
$X\in\Sigma(k, g, S)$ are defined over $k(g, S)$.

Pick $X\in\Sigma(k, g, S)$.  As $2$ is a unit in $\oo_S$ and $X/k$ has good
reduction outside $S$, for every maximal ideal $\p$ of $\oo_S$ the reduction
mod $\p$ map induces a bijection between the Weierstrass points of
$
X/\overline{k}
$
and those of
$
X/\overline{k_\p}
$.
Then the (conjugacy class of the) inertia group $I_\p$ of $\p$ in 
$
\text{Gal}(\overline{k}/k)
$
acts trivially on the Weierstrass points of $X/\overline{k}$, since $I_\p$ 
acts trivially on $k_\p$.  There are exactly $2g+2$ such points, so we are
done.
\end{proof}

Let $g\ge 2$. We now show that $\Sigma_0(k, g, S)$ is finite for all $k$ and
for all $S$ satisfying (i) and (ii).
Assume as inductive hypothesis that this is so for curves of positive genus
$<g$, the case $g=2$ being taken care of by Shafarevich
\cite{shaf}. Let $X/k$ be a pointed
hyperelliptic
curve of genus $g$.  By
\cite[Cor.~2.10]{lock} and by our choice of $S$, there exists a model for
$X$ of the form
$
\mathcal C_X: y^2 + y Q(x) = P(x)
$
with
$
Q, P\in\oo_S[x], \deg P = 2g+1
$,
and
$
\deg Q\le g
$,
such that for any maximal ideal $\p\subset\oo_k$ not in $S$,  the curve $X$
has bad reduction at $\p$ if and only if the
discriminant
$
\Delta(\mathcal C_X) := 2^{4g}\text{disc}(P(x) + \frac{1}{4}Q(x)^2)
$
of the model $\mathcal C_X$ (cf.~\cite[$\S 1$]{lock}) is divisible by $\p$,
viewed as a maximal ideal in $\oo_S$.
Furthermore, since $2$
is a unit in $\oo_S$, by completing the square we can assume that $Q=0$.
So if $X$ has good reduction outside $S$, then the principal ideal
$\text{disc}(P)\oo_S$ is equal to $\oo_S$,
and the polynomial $P\in k[x]$ splits completely over
$K' = k(g, S)$.
Without loss of generality we can assume that $P(0)\not=0$.  Let
$
R(z) = z^{2g+2}P(1/z)
$.
Then $y^2=R(z)$ is another model of $X$ over $\oo_S$, we have $\deg R=2g+2$,
and $R$ splits completely over $K'$.  Thus 
we can write (not uniquely) $R=R_1 R_2$, where
$R_i\in K'[x]$ and $\deg R_1 = 3, \deg R_2 = 2g-1$. 
Denote by $S'$ the set of primes of $K'$ lying above those in $S$. From the
discriminants of $R_1$ and $R_2$ we see that the two
$K'$-rational hyperelliptic curves
$
X_i: y^2 = R_i(z)
$
have good reduction outside $S'$.   Moreover, both $X_i$ are
pointed over $K'$, and their genus are both positive and
are less than that of $C$.  By induction there are at most
finitely many  $K'$-isomorphism classes of such $X_i/K'$.

Denote by $X_3/K'$ the fiber product
$
X_1\times_{\mathbf P^1_{K'}} X_2
$
with respect to the hyperelliptic involutions
$
X_i\rarr\mathbf P^1_{K'}
$.
An affine coordinate ring for $X_3/K'$ is
$
A = K'[z, \sqrt{R_1}]\otimes_{K'[z]}K'[z,\sqrt{R_2}]
$.
It contains a subring isomorphic to 
$
K'[z,\sqrt{R_1 R_2}]
$,
which is an affine coordinate ring for $X/K'$.  Moreover, $K'$ is
algebraically closed in the field of fraction of $A$.  Consequently, if we
denote by $X'_3/K'$ the
smooth model of $X_3/K'$, there is a non-trivial $K'$-morphism
$
X'_3\rarr X
$.
The Riemann-Hurwitz formula then implies that $X'_3$
has genus $\ge 2$ since $X$ does.

To recapitulate, with $k, S$ and $g\ge 2$ as above, there exists a finite
collection
$\Xi(K', g, s)$ of smooth curves  over $K'$ of genus $\ge 2$, such that every
$
X\in\Sigma_0(k, g, S)
$
when based changed to $K'$ is the $K'$-image of some $X'_3\in\Xi(K', g, S)$.
By the Theorem of de Franchis \cite{de}, there are at most finitely many
$\overline{K'}$-isomorphism classes of
pairs $(\pi, Z)$, where $Z$ is a smooth projective curve of genus $\ge 2$ over
$\overline{K'}$, and $\pi: X_3'\rarr Z$ is a non-trivial map over
$\overline{K'}$.  The same
finiteness conclusion then
holds if we replace $\overline{K'}$ by $K'$, since the automorphism group of
any curve of genus $\ge 2$ is finite.  Consequently,
$\Sigma_0(k, g, S)$ falls into finitely many $K'$-isomorphism classes.
Invoke the $H^1$-finiteness argument in the proof of the Lemma and we get the
finiteness of $\Sigma_0(k, g, S)$ for all $k$ and $S$.  Apply the Lemma and we
get the finiteness of $\Sigma(k, g, S)$ for all $k$ and $S$, as desired.

\bibliographystyle{amsalpha}

\end{document}